\documentclass[12pt,twoside,doublespace]{article}

\usepackage{amssymb}
\usepackage{graphicx}

\def\smskip{\par\vskip 5 pt}
\def\QED{\hfill $\Box$\smskip}
\newtheorem{theorem}{Theorem}
\newtheorem{lemma}{Lemma}
\newtheorem{proposition}{Proposition}

\begin{document}

\begin{center}

\vspace{35pt}

{\Large \bf A Simple Adaptive Step-size Choice }

\vspace{5pt}

{\Large \bf for Iterative Optimization Methods}

\vspace{35pt}

{\sc I.V.~Konnov\footnote{\normalsize E-mail: konn-igor@ya.ru}}

\vspace{35pt}

{\em  Department of System Analysis
and Information Technologies, \\ Kazan Federal University, ul.
Kremlevskaya, 18, Kazan 420008, Russia.}

\end{center}

\vspace{35pt}

\begin{abstract}
We suggest a simple adaptive step-size procedure, which does not require any line-search,
for a general class of nonlinear optimization methods
and prove convergence of a general method under mild assumptions.
In particular, the goal function may be non-smooth and non-convex.
Unlike the descent line-search methods, it does not require monotone decrease of the goal function 
values along the iteration points and 
reduces the implementation cost of each iteration essentially.
The key element of this procedure consists in
inserting a majorant step-size sequence such that the next element is taken only
if the current iterate does not give a sufficient descent.
Its applications yield in particular
a new gradient projection method for smooth constrained optimization problems
and a new projection type method for minimization of the gap function of a
general variational inequality.
Preliminary results of computational experiments confirm efficiency of
the proposed modification.


\end{abstract}


\newpage

%1111111111111111111111111111111111111111111111111111111111111111111111111

\section{Introduction} \label{sc:1}

Iterative methods are utilized for various difficult problems
whose solution in the closed form either is not known or
has significant computational drawbacks. For instance,
if a system of linear equations has large dimensionality and inexact data,
it is better to apply one of well known iterative methods.
These methods are a standard tool for nonlinear constrained
optimization problems; see e.g. \cite{Pol87,Ber99}.
One of the most popular approaches to generation of an iterative sequence
consists of solution of direction finding
and step-size choice subproblems at each iteration.

During rather long time, most efforts
were concentrated on developing more powerful and rapidly convergent methods,
such as Newton and interior point type ones, which admit
complex transformations at each iteration and attain high accuracy
of approximations. That is, the direction finding subproblem was considered
as the main one, whereas the step-size was chosen by one of the few well known
procedures; see e.g. \cite{Ber99,BV09,Pat99}.
 However, new significant areas of applications related to
data processing in information and communication systems,
having large dimensionality and inexact data together with
scattered necessary information force one to avoid
complex transformations and apply mostly simple
methods such as the projection and linearization type methods,
whose iteration computation expenses and accuracy requirements
are rather low; see e.g. \cite{Bur98,BT09,CBS14,SFSPP14}.
 Therefore, one is also interested in suggesting new
step-size choice rules that reduce the total computational expenses
of the method.

In fact, the existing rules are not completely satisfactory.
The exact or approximate one-dimensional
minimization line-search requires significant computational expenses per iteration
especially in the case where calculation of the function value is
almost similar to calculation of its derivative (gradient) and needs
solution of complex auxiliary problems; see e.g. \cite{Pol87,BV09}.
In order to remove the  line-search, one can calculate the step-size value via
utilization of a priori information such as Lipschitz type constants for the gradient, but then one must
take only some their inexact estimates, which leads to slow convergence.
This is also the case for the known divergent series rule; see e.g. \cite{Pol87,BT09}.

In this paper, we suggest a new simple adaptive step-size procedure for a general class of
iterative optimization methods, which does not require any line-search.
In particular, this procedure can be applied to the known projective optimization methods.
In creation of the adaptive step-size rule we follow the approach from \cite{Kon18} where
a step-size procedure for the conditional gradient method was proposed.
However, the procedure in \cite{Kon18} admits only decrease of the step-size
and can not be extended to the other optimization methods since it requires
the boundedness of the feasible set and is adjusted to the set-valued solution mapping of
the direction finding subproblem. Our new step-size procedure admits different changes of the step-size
and wide variety of implementation rules.
The key element of this procedure consists in
inserting a majorant step-size sequence tending to zero.
However, we do not change these majorant values continuously
such that we take the next element only
if the current iterate has not given a sufficient descent.
In such a way, the procedure takes into account behavior of the iteration sequence.
We show that this strategy can be implemented within a rather general framework
of iterative solution methods applied both
for smooth and non-smooth optimization problems and variational inequalities.
It does not utilize a priori information such as Lipschitz constants
of the gradient, besides, the Lipschitz continuity of the
gradient of the goal function is not necessary for convergence of the method.
Preliminary results of computational experiments confirm efficiency of
the proposed modification.

The remainder of the paper is organized as follows. Section \ref{sc:2} contains
necessary definitions and properties from the theory of
non-smooth optimization and set-valued analysis. In Section \ref{sc:3},
we describe a general framework of iterative methods with
the new step-size procedure applied to
constrained non-smooth and non-convex optimization problems
and prove its convergence.
Afterwards, we describe its specializations as projection type methods.
Namely, a new gradient projection method for smooth constrained optimization problems is given
in Section \ref{sc:4}, whereas a new projection type method for a
variational inequality with a general non-integrable mapping, which is re-formulated as
a constrained gap function minimization problem, is given
in Section \ref{sc:5}. We show that they both fall into the
general framework of Section \ref{sc:3} and obtain their convergence
directly from the basic convergence theorem of Section \ref{sc:3}.
Section \ref{sc:6} contains
results of preliminary computational experiments.

%22222222222222222222222222222222222222222222222222222222222222222222222222222

\section{Basic Preliminaries}\label{sc:2}

We intend to develop the new method for a wide class of optimization problems
whose goal functions can be non-smooth and non-convex.
For this reason, we first recall some concepts and properties from Non-smooth Analysis; see \cite{Cla83}
for more details. If some function $f : \mathbb{R}^{n} \to \mathbb{R}$
is Lipschitz continuous in a neighborhood of
any point $x$ of a set $X$, it is called locally Lipschitz on $X$.
Then we can define its generalized gradient set at $x$:
$$
\partial^{\uparrow} f(x)=\{ g \in \mathbb{R}^{n} \ : \
\langle g , p \rangle \le f^{\uparrow}(x;p)   \quad \forall p \in
\mathbb{R}^{n}\},
$$
which must be non-empty, convex and closed. Here $f^{\uparrow}(x;p)$ denotes the upper Clarke
derivative:
$$
  f ^{\uparrow}(x;p) = \limsup _{y \to x, \alpha \searrow 0}
 ((f (y+ \alpha p) - f (y))/ \alpha ).
$$
It follows that
$$
  f ^{\uparrow}(x,p) = \sup _{g \in \partial^{\uparrow} f(x)}
            \langle g, p \rangle.
$$
In general, the locally Lipschitz function $f$ need not be differentiable.
At the same time,  it has the gradient $\nabla f(x)$ a.e. in $X$,
furthermore, it holds that
\begin{equation} \label{eq:2.1}
  \partial^{\uparrow} f(x)
={\rm conv} \left\{\lim \limits_{y\rightarrow x}\nabla f(y)\ : \
y\in X_{f}, \ y\notin S \right\},
\end{equation}
where $X_{f}$ denotes  the set of points where $f$ is
differentiable, and $S$ denotes  an arbitrary subset of measure
zero. If $f$ is convex, then  $\partial^{\uparrow} f(x)$ coincides
with the subdifferential $\partial f(x)$ in the sense of Convex Analysis,
i.e.,
$$
 \partial f(x) = \{  g \in \mathbb{R}^{n}  \ : \ f(y)-
  f(x) \geq  \langle g, y-x \rangle   \quad \forall y \in \mathbb{R}^{n} \}.
$$
In this case the upper
derivative coincides with the usual directional derivative:
\begin{equation} \label{eq:2.2}
  f ^{\uparrow}(x;p) = f'(x;p).
\end{equation}
Also, if $f$ is differentiable at $x$, (\ref{eq:2.2}) obviously
holds and we have
$$
  f'(x;p) = \langle  \nabla f(x), p \rangle  \
   {\rm and } \  \partial ^{\uparrow} f(x) = \{ \nabla f(x) \};
$$
cf. (\ref{eq:2.1}).

The general optimization problem consists in finding
the minimal value of some goal function $f : \mathbb{R}^{n} \to \mathbb{R}$
on a feasible set $D \subseteq \mathbb{R}^{n}$. For brevity, we
write this problem as
\begin{equation} \label{eq:2.3}
 \min \limits _{x \in D} \to f(x),
\end{equation}
its solution set is denoted by $D^{*}$ and the optimal value of the
function by $f^{*}$, i.e. $f^{*} = \inf\limits_{x \in D} f(x)$.
We will use the following first set of basic
 assumptions for problem (\ref{eq:2.3}).

{\bf (A1)}  {\em The set  $D$ is nonempty, convex, and closed,
the function $f : \mathbb{R}^{n} \to \mathbb{R}$ is locally Lipschitz on $D$.}

{\bf (A2)} {\em There exists a number $\gamma > f^{*}$ such that the set
$$
D_{\gamma}=\left\{ x \in D \ : \ f (x) \leq \gamma \right\}
$$
is bounded.}

Clearly, {\bf (A2)} is a general coercivity condition that
is necessary in the case where the set $D$ is unbounded.
If {\bf (A1)} and {\bf (A2)} hold, problem (\ref{eq:2.3}) has a solution.
This means that we intend to present a method for
non-smooth and non-convex optimization problems.

Together with problem (\ref{eq:2.3}) we will consider
the following set-valued variational inequality (VI for short): Find a point
$x^{*} \in D$ such that
\begin{equation} \label{eq:2.4}
\exists g^{*} \in \partial^{\uparrow} f (x^{*}),
\quad  \langle g^{*},x-x^{*} \rangle \geq 0 \quad \forall x \in D.
\end{equation}
We denote by $D^{0}$ the solution set of VI (\ref{eq:2.4}).
Solutions of VI (\ref{eq:2.4}) are called {\em stationary points}
of (\ref{eq:2.3}) due to the known necessary optimality condition;
see e.g. \cite{Cla83,Kon13}.

%============================pro:2.1======================================

\begin{proposition} \label{pro:2.1} {\em Let {\bf (A1)} hold.
Then each solution of problem (\ref{eq:2.3}) is a solution of VI (\ref{eq:2.4}).}
\end{proposition}

The reverse implication needs additional conditions.
We recall that a function $\varphi : \mathbb{R}^{n}
\to \mathbb{R}$ is called

(a) {\em pseudo-convex} on a set $X$, if
for each pair of points $x, y \in X$, we have
$$
    \varphi' (x; y - x ) \geq 0 \
      \Longrightarrow \
    \varphi (y) \geq \varphi (x);
$$

(b) {\em semi-convex} (or {\em upper pseudo-convex}) if
for each pair of points $x, y \in X$, we have
$$
    \varphi^{\uparrow}(x; y - x ) \geq 0 \
      \Longrightarrow \
    \varphi (y) \geq \varphi (x);
$$
see \cite{Mif77} and also \cite{Kon13}. In case (\ref{eq:2.2}), these concepts coincide, but
in general (b) implies (a). Besides, the class of convex functions is strictly contained
in that of pseudo-convex functions. If {\bf (A1)} holds and $f$ is semi-convex, then  each solution of VI
(\ref{eq:2.4}) clearly solves  problem  (\ref{eq:2.3}), i.e. $D^{*}=D^{0}$.

We need several continuity properties of set-valued mappings; see
e.g. \cite{Nik68,Kon13}. Here and below $\Pi(A)$ denotes  the family
of all nonempty subsets of a set $A$.

Let $X$ be a convex set in $\mathbb{R}^{n}$.
A set-valued mapping $Q: X \to \Pi (\mathbb{R}^{n})$ is said to be

 (a) {\em upper semicontinuous (u.s.c.)},
 if for each point $y \in X$ and for each open set $U$ such
that $Q(y) \subset U$, there is a neighborhood $Y$ of $y$ such that
$Q(z) \subset U$ whenever $z \in X \cap Y$;

 (b) {\em closed}, if for each pair of sequences
 $\{ x^{k} \} \to x$,   $\{ q^{k} \} \to q$
such that $x^{k} \in X$ and $q^{k} \in Q(x^{k})$, we have $q \in
Q(x)$;

 (c) a {\em K-mapping} (Kakutani-mapping), if it is
 u.s.c. and has nonempty, convex, and compact values.

It is known (see e.g. \cite[Chapter 1, Lemma 4.4]{Nik68}), that each
 u.s.c. mapping with closed values is closed and that each
 closed mapping which maps any compact set into a  compact set is  u.s.c.
Also, if a function $f:Y \rightarrow \mathbb{R}$ is locally Lipschitz
on an open convex set $Y$, then
$\partial^{\uparrow} f$ is a K-mapping on $Y$; see
\cite[Section 2.1]{Cla83}.

We shall also  use the mean value theorem by G.~Lebourg for locally
Lipschitz functions.

%1===========================pro 2.2====================================================

\begin{proposition} \label{pro:2.2} {\em
\cite[ Theorem 2.3.7]{Cla83} Let $x$ and $y$ be given points in
$\mathbb{R}^{n}$ and let $f:\mathbb{R}^{n}\rightarrow \mathbb{R}$ be
a Lipschitz continuous function on an open set containing the segment $[x,
y]$. Then there exists a point $z\in (x, y)$ such that
$$
f(y)-f(x)\in \langle \partial^{\uparrow} f(z), y-x \rangle .
$$}
\end{proposition}

%33333333333333333333333333333333333333333333333333333333333333333333333333

\section{The Basic Method and Its Convergence} \label{sc:3}

We first describe conditions for the solution mapping of
the direction finding subproblem.

{\bf (A3)} {\em There exists a single-valued mapping $x \mapsto y(x)$,
which maps the set $D$ into $D $ such that

(i)  it is continuous;

(ii)  $\bar x= y(\bar x) $ if and only if  $\bar x $ is a solution of VI (\ref{eq:2.4});

(iii)  for each $x\in D$ and for all $g\in \partial^{\uparrow} f(x)$ it holds that
\begin{equation} \label{eq:3.1}
 \langle g, y(x)-x \rangle \leq -\tau \|y(x)-x\|^{2}
\end{equation}
for some $\tau>0 $.
}

We now describe the basic method  for problem (\ref{eq:2.3}), which involves a
simple adaptive step-size procedure without line-search.

%===========================Method (SBM)==================================
\medskip
\noindent {\bf Method (SBM).}

{\em Step 0:} Choose a point $x^{0}\in D_{\gamma}$, a number $\beta \in (0,1)$ and a sequence
$\{\tau _{l}\} \to 0$, $\tau _{l} \in (0, 1)$. Set $k=0$, $l=0$, $u^{0}=x^{0}$,
choose a number $\lambda_{0} \in (0, \tau _{0}]$.

{\em Step 1:} Take a point $y^{k}=y(x^k)$.
If $y^{k}=x^{k}$, stop. Otherwise set $d^{k}=y^{k}-x^{k}$
and $z^{k+1}=x^{k}+\lambda_{k}d^{k}$.

{\em Step 2:} If
\begin{equation} \label{eq:3.2a}
  f(z^{k+1}) \leq f(x^{k})-\beta \lambda_{k}\|d^{k}\|^{2},
\end{equation}
take $\lambda_{k+1} \in [\lambda_{k}, \tau _{l}]$, set $x^{k+1}=z^{k+1}$
and go to Step 4.

{\em Step 3:} Set $\lambda'_{k+1} = \min \{\lambda_{k}, \tau _{l+1}\}$,
$l=l+1$ and take $\lambda_{k+1} \in (0,\lambda'_{k+1}]$.
If $f (z^{k+1}) \leq \gamma$, set $x^{k+1}=z^{k+1}$ and go to Step 4.
Otherwise set $x^{k+1}=u^{k}$, $u^{k+1}=u^{k}$, $k=k+1$ and go to Step 1.

{\em Step 4:} If $f(x^{k+1}) < f(u^{k})$, set $u^{k+1}=x^{k+1}$,
  $k=k+1$ and go to Step 1.
\medskip

Therefore, (SBM) in fact represents a general framework of iterative methods for
nonlinear optimization problems.
Observe that the sequence $\{u^{k}\}$ simply contains the best current points of
the sequence $\{x^{k}\}$, i.e.
$$
f(u^{k})=\min_{0 \leq i \leq k} f(x^{i}).
$$
Due to {\bf (A3)}, termination of (SBM) yields a point of $D^{0}$. Hence,
we will consider only the case where the sequence $\{x^{k}\}$ is infinite.

%3=========================thm 3.1============================================

\begin{theorem} \label{thm:3.1} {\em
Let the assumptions {\bf (A1)}--{\bf (A3)} be fulfilled and $\beta< \tau$. Then:

(i) The sequence $\{x^k \}$ has a limit point, which belongs to the set
$D^{0}$.

(ii)  If $D^{*}=D^{0}$, then all the limit points of the sequence $\{x^k \}$
belong to the set $D^{*}$, besides, we have
\begin{equation}\label{eq:3.2}
\lim \limits_{k\rightarrow \infty} f(x^{k})=f^{*}.
\end{equation}}
\end{theorem}
{\bf Proof.}  First we observe that the sequence $\{x^{k}\}$ belongs to the bounded
set $D_{\gamma}$ and must have limit points. By {\bf (A3)}, so is
the sequence $\{y^{k}\}$, hence $\{d^{k}\}$.
Next, we take the subsequence of indices $\{i_{s}\}$
 such that
 \begin{eqnarray}
  & &   f (z^{i_{s}+1}) > \gamma, \ f (x^{i_{s}}) \leq \gamma, \label{eq:3.3}\\
   &&  f(z^{i_{s}+1}) > f(x^{i_{s}})-\beta \lambda_{i_{s}}\| d^{i_{s}}\|^{2}, \
   z^{i_{s}+1}=x^{i_{s}}+ \lambda_{i_{s}} d^{i_{s}}.
    \label{eq:3.4}
\end{eqnarray}

Let us consider several possible cases.

\textit{Case 1: The subsequence $\{x^{i_{s}}\}$ is infinite.} \\
Take  an arbitrary limit point $x'$  of the subsequence
$\{x^{i_{s}}\}$.  Without loss of generality we can suppose that
$$
\lim \limits_{s\rightarrow \infty }x^{i_{s}}= x' \ \mbox{and} \ \lim
\limits_{s\rightarrow \infty }y^{i_{s}}= y',
$$
where $y'=y( x')$ by {\bf (A3)}. Note that
$$
 \lambda_{i_{s}} \in (0,\tau_{l_{s}}], \ \lambda_{i_{s}+1} \in (0,\tau_{l_{s}+1}],
$$
for some infinite subsequence of indices $\{l_{s}\}$ where
$$
\lim \limits_{s\rightarrow \infty}\tau_{l_{s}}=0.
$$
Since the sequence $\{d^{i_{s}}\}$ is bounded,
 the limit points of the subsequences $\{x^{i_{s}}\}$ and
$\{z^{i_{s}+1}\}$ coincide due to (\ref{eq:3.4}). From (\ref{eq:3.3})
we now obtain
\begin{equation}\label{eq:3.5}
 f(x')=\gamma>f^{*}.
\end{equation}
Applying Proposition \ref{pro:2.2} in (\ref{eq:3.4}), we have
$$
\langle g^{i_{s}}, d^{i_{s}} \rangle \geq -\beta \|d^{i_{s}}\|^{2}
$$
for some  $g^{i_{s}}\in \partial^{\uparrow} f(x^{i_{s}}+\theta _{i_{s}}\lambda _{i_{s}}d^{i_{s}})$,
and $\theta _{i_{s}}\in (0,1)$. Taking the limit $s\rightarrow \infty $  gives
$$
\langle g', y'-x' \rangle \geq -\beta \|y'-x'\|^{2}
$$
for some  $g'\in \partial^{\uparrow} f(x')$.
Using (\ref{eq:3.1}), we obtain
$$
\beta \| y'-x'\|^{2}\geq \tau \|y'-x'\|^{2},
$$
i.e. $y(x')=x'$, hence
\begin{equation}\label{eq:3.6}
 x' \in D^{0}.
\end{equation}
Therefore, assertion (i) is true in this case.

\textit{Case 2: The subsequence $\{x^{i_{s}}\}$ is finite.} \\
Without loss of generality we can then suppose that
$z^{k}=x^{k}$ for each number $k$. The further proof depends on the properties of
the sequence  $\{\lambda _{k}\}$.

\textit{Case 2a: The number of changes of the index $l$ is finite.} \\
Then we have $\lambda _{k} \geq \bar \lambda>0$ for $k$ large enough, hence
(\ref{eq:3.2a}) gives
$$
 f(x^{k+1}) \leq f(x^{k})-\beta \lambda_{k}\|d^{k}\|^{2}
            \leq f(x^{k})-\beta \bar \lambda \|d^{k}\|^{2}
$$
for $k$ large enough. Since $f(x^{k}) \geq f^{*}> -\infty$, we must have
\begin{equation} \label{eq:3.7}
\lim \limits_{k\rightarrow \infty }f(x^{k})=\mu
\end{equation}
and
\begin{equation} \label{eq:3.8}
\lim \limits_{k\rightarrow \infty}\|y^{k}-x^{k}\|=0.
\end{equation}
Let $x''$ be an arbitrary limit point of the sequence  $\{x^{k}\}$.
From (\ref{eq:3.8}) and {\bf (A3)} we now have
$$
y(x'')=x'',
$$
which gives  $x'' \in D^{0}$.  Hence
in this case all the limit points of the sequence $\{x^k \}$
belong to the set $D^{0}$. Therefore,
assertion (i) is true in this case.

\textit{Case 2b: The number of changes of the index $l$ is infinite.} \\
Then there exists an infinite subsequence of indices $\{k_{l}\}$
 such that
\begin{equation} \label{eq:3.9}
f (x^{k_{l}}+\lambda _{k_{l}}d^{k_{l}})-f (x^{k_{l}})=f
(x^{k_{l}+1})-f (x^{k_{l}}) > -\beta \lambda _{k_{l}}\| d^{k_{l}}\|^{2}.
\end{equation}
besides,
$$
 \lambda_{k_{l}} \in (0,\tau_{l}], \ \lambda_{k_{l}+1} \in (0,\tau_{l+1}],
$$
and
$$
\lim \limits_{l\rightarrow \infty}\tau_{l}=0.
$$
Let $\bar x$ be an arbitrary limit point  of this subsequence
$\{x^{k_{l}}\}$.  Without loss of generality we can suppose that
$$
\lim \limits_{l\rightarrow \infty }x^{k_{l}}=\bar x \ \mbox{and} \ \lim
\limits_{l\rightarrow \infty }y^{k_{l}}=\bar y.
$$
where $\bar y=y(\bar x)$. Applying Proposition \ref{pro:2.2} in (\ref{eq:3.9}), we have
$$
\langle g^{k_{l}}, d^{k_{l}} \rangle \geq -\beta \|d^{k_{l}}\|^{2}
$$
for some  $g^{k_{l}}\in \partial^{\uparrow} f(x^{k_{l}}+\theta _{k_{l}}\lambda _{k_{l}}d^{k_{l}})$,
and $\theta _{k_{l}}\in (0,1)$. Since $ \lambda _{k_{l}} \to 0$ as $l\rightarrow \infty$,
taking the limit $l\rightarrow +\infty $ gives
$$
\langle \bar  g, \bar y-\bar x \rangle \geq -\beta \|\bar y-\bar x\|^{2}
$$
for some  $\bar g\in \partial^{\uparrow} f(\bar x)$.
Using (\ref{eq:3.1}), we obtain
$$
\beta \| \bar y-\bar x\|^{2}\geq \tau \|\bar y-\bar x\|^{2},
$$
i.e. $\bar x=y(\bar x)$, hence
$ \bar x\in D^{0}$. Therefore, all the limit
points of the subsequence $\{x^{k_{l}}\}$ belong to the set $D^{0}$.
Since $x^{k_{l}+1}=x^{k_{l}}+\lambda _{k_{l}}d^{k_{l}}$, $ \lambda
_{k_{l}} \to 0$, and the sequence  $\{d^{k_{l}}\}$ is bounded,  the
limit points of the subsequences $\{x^{k_{l}}\}$ and
$\{x^{k_{l}+1}\}$ coincide and  all they belong to the set $D^{0}$.
We conclude that assertion (i) is also true in this case.

We now suppose in addition that $D^{0}=D^{*}$. Then relations
(\ref{eq:3.5}) and (\ref{eq:3.6}) become inconsistent, hence
Case 1 is impossible. This means that the subsequence $\{x^{i_{s}}\}$
is always finite. In Case 2a we now have $\mu=f^{*}$
in (\ref{eq:3.7}), which gives (\ref{eq:3.2}).
We conclude that assertion (ii) holds true in this case.

In Case 2b,  the limit points of the subsequences $\{x^{k_{l}}\}$ and
$\{x^{k_{l}+1}\}$ coincide and all they now belong to the set $D^{*}$.
For any index $k$ we define the index $m(k)$ as follows:
$$
 \ m(k)= \max\{j \ : \ j \leq k, \ f (x^{j})-f (x^{j-1}) > -\beta
\lambda _{j-1}\| d^{j-1}\|^{2}\},
$$
i.\,e. $m(k)$ is the closest to $k$ but not greater index from the
subsequence $\{x^{k_{l}+1}\}$. This means that $m(k)=k$ if $f (x^{k})-f
(x^{k-1}) > -\beta \lambda _{k-1}\| d^{k-1}\|^{2}$.
By definition, we have
\begin{equation} \label{eq:3.10}
f(x^{k})\leq f (x^{m(k)}).
\end{equation}
Let now $x^{*}$ be an arbitrary
limit point  of the sequence $\{x^{k}\}$, i.e. $\lim
\limits_{s\rightarrow \infty }x^{t_{s}}=x^{*}$. Create the corresponding infinite
subsequence $\{x^{m(t_{s})}\}$. From (\ref{eq:3.10}) we have $f^{*} \leq f
(x^{t_{s}})\leq f (x^{m(t_{s})})$, but all the limit points of the
sequence $\{x^{m(t_{s})}\}$ belong to the set $D^{*}$ since it is
contained in the sequence $\{x^{k_{l}+1}\}$. Choose any limit point
$\tilde x$ of $\{x^{m(t_{s})}\}$.  Then, taking a subsequence if necessary we
obtain
$$
f^{*} \leq f (x^{*})\leq f (\tilde x) = f^{*}.
$$
therefore $x^{*} \in D^{*}$. This means that all the limit points of the
sequence $\{x^{k}\}$ belong to the set  $D^{*}$ and that
(\ref{eq:3.2}) holds true. We conclude that
assertion (ii) is  also true.
\QED 

The method can be simplified in the case where the set $D$ is bounded.
Then we can set $\gamma=+ \infty$ and remove
all the calculations of the sequence $\{u^{k}\}$.
It is easy to verify that all the assertions
 of Theorem \ref{thm:3.1} remain true.

%444444444444444444444444444444444444444444444444444444444444444444444444444444444444444444

\section{Application to Smooth Optimization Problems}\label{sc:4}

From the results of Section \ref{sc:3} it follows that we can
create a number of new solution methods for optimization problems.
It suffices to take an optimization problem that satisfies conditions
{\bf (A1)} and {\bf (A2)} and a method whose solution mapping of
the direction finding subproblem satisfies condition
{\bf (A3)}. Then we place this method in the framework of
(SBM) and obtain its convergence properties
directly from Theorem \ref{thm:3.1}.

We illustrate diversity of possible specializations
of (SBM) by only two basic examples. In this section,
we take the well known class of smooth constrained
optimization problems.

{\bf (A1${}'$)}  {\em The set  $D$ is nonempty, convex, and closed,
the function $f : \mathbb{R}^{n} \to \mathbb{R}$ is continuously differentiable on $D$.}

Clearly, {\bf (A1${}'$)}  implies {\bf (A1)}. Let $\pi _{X} (x)$ denotes the
projection of $x$ onto a set $X$. Fix a number $\alpha >0$ and define
the mapping $y_{\alpha }(x)=\pi_{D}[x-\alpha^{-1} \nabla f(x)]$ on the set $D$.
Then setting $y(x)=y_{\alpha }(x)$ in (SBM), we obtain a new version of
the gradient projection method for problem (\ref{eq:2.3}). We call it (GPMS) for brevity.

We now utilize the well known properties of mapping $x \mapsto y_{\alpha }(x)$;
see e.g. \cite{Pat98} and \cite[Lemma 9.5]{Kon13}.

%4=========================lm 4.1====================================================
\begin{lemma} \label{lm:4.1} {\em Let the assumptions in {\bf (A1${}'$)} be fulfilled. Then:

(a) $\bar x= y_{\alpha }(\bar x) $ if and only if  $\bar x \in D^{0}$;

(b)  The mapping $x\mapsto y_{\alpha }(x)$ is continuous on $D$;

(c) For any point $x\in D$ it holds that
$$
\langle  \nabla f(x), y_{\alpha }(x)-x \rangle   \leq - \alpha \|y_{\alpha }(x)-x\|^{2}.
$$}
\end{lemma}
Therefore, the assumptions in {\bf (A3)} are fulfilled with $\tau=\alpha$ and we can
obtain the convergence result for (GPMS) directly from Theorem \ref{thm:3.1}.

%4=========================thm 4.1============================================

\begin{theorem} \label{thm:4.1}{\em
Let assumptions {\bf (A1${}'$)} and {\bf (A2)} be fulfilled.
If we apply (GPMS) with $\beta< \alpha$, then the following assertions are true:

 (i) The termination of (GPMS) yields a point of $D^{0}$.

(ii) The sequence $\{x^k \}$ has a limit point, which belongs to the set
$D^{0}$.

(iii)  If $D^{*}=D^{0}$, then all the limit points of the sequence $\{x^k \}$
belong to the set $D^{*}$ and (\ref{eq:3.2}) holds.}
\end{theorem}

Observe that the choice $\alpha \geq  1$ allows us to take an arbitrary value
$\beta \in (0,1)$ for convergence. The above method can be extended to the case where
$$
f(x)=\mu(x)+\eta (x),
$$
where $\mu$ is continuously differentiable and $\eta$ is convex, but non-differentiable.
We have to replace the projection mapping with the proximal mapping with respect to $\eta$
and apply the so-called splitting or proximal gradient iteration.
A number of these splitting based descent methods were proposed for such composite non-smooth
optimization problems; see e.g. \cite{Pat98,Pat99,BT09,Kon13,CBS14}.
Similarly, we can create new versions of splitting methods
if we place the corresponding splitting direction finding mapping in the (SBM) framework.

%555555555555555555555555555555555555555555555555555555555555555555555555555

\section{Application to Non-smooth Variational Inequality Problems}\label{sc:5}

In this section, we take a variational inequality
whose underlying mapping is (strongly) monotone, but non-integrable and non-smooth
in general. Given a convex set $D$  in $\mathbb{R}^{n}$ and  a single-valued
mapping $G : D \rightarrow \mathbb{R}^{n}$,
one can define the custom {\em variational inequality} (VI for
short): Find $x^{*} \in D$ such that
\begin{equation} \label{eq:5.1}
 \langle G(x^{*}), x-x^{*} \rangle \geq 0 \quad \forall x\in D.
\end{equation}
This problem has a great number of applications in different fields,
the theory and methods for VIs are investigated very
extensively and great advances were made in this field;
e.g. see \cite{Pat99,Kon01,FP03} and the references therein.
We denote by $D^{e}$ the solution set of VI (\ref{eq:5.1})
and  consider this problem  under the
following basic assumptions.

{\bf (A1${}''$)} \ {\em $D$ is a nonempty, closed, and convex set in
$\mathbb{R}^{n}$, the mapping $G : Y \rightarrow \mathbb{R}^{n}$
satisfies the Lipschitz condition in a neighborhood of each point of
an open convex set $Y$ such that $D\subset Y$.}

Fix a number $\alpha >0$ and define the usual gap function
\begin{equation} \label{eq:5.2}
  \varphi _{\alpha }(x) = \max _{y \in D}  \{\langle G(x), x - y \rangle -
     0.5 \alpha  \| x- y \|^{2}\};
\end{equation}
e.g. see \cite{Pat99}.
Then  there exists a unique element $ y_{\alpha }(x) \in D$
such that
$$
 \varphi _{\alpha }(x) = \langle G(x), x - y_{\alpha }(x) \rangle -
     0.5 \alpha  \| x- y_{\alpha }(x) \|^{2},
$$
moreover, $ y_{\alpha }(x) =\pi_{D}[x-\alpha^{-1} G(x)]$. Thus, we again have
a single-valued mapping $x \mapsto y_{\alpha }(x)$ on $D$.
Then, setting $y(x)=y_{\alpha }(x)$ and $f= \varphi _{\alpha }$ in (SBM), we obtain a new version of
the projective gap function method for problem (\ref{eq:5.1}). We call it (GFPMS) for brevity.
We observe that the descent method with Armijo line-search  was proposed for
this non-smooth VI in \cite{Kon06f}.

First of all we replace VI (\ref{eq:5.1}) with
the optimization problem
\begin{equation} \label{eq:5.3}
\begin{array}{c}
\displaystyle \min \limits_{x\in D}\rightarrow \varphi _{\alpha}(x).
\end{array}
\end{equation}

From the definition of the function $ \varphi _{\alpha }$ in (\ref{eq:5.2}) we can
easily deduce that it is always nonnegative and that the optimal value in (\ref{eq:5.3}) is zero.
We shall utilize the other known properties of the mapping $x \mapsto y_{\alpha }(x)$;
see  \cite{Kon06f,Kon12}.

%5=========================lm 5.1====================================================
\begin{lemma} \label{lm:5.1} {\em Let the assumptions in {\bf (A1${}''$)} be fulfilled. Then:

(a)  VI (\ref{eq:5.1}) is equivalent to  problem (\ref{eq:5.3});

(b)  $\bar x= y_{\alpha }(\bar x) $ if and only if  $\bar x \in D^{e}$;

(c) The mapping $x\mapsto y_{\alpha }(x)$ is continuous on $D$.}
\end{lemma}

Therefore, the optimization problem (\ref{eq:5.3}) is an equivalent re-formulation of
VI (\ref{eq:5.1}). Although the  gap function $\varphi _{\alpha}$ is non-smooth and non-convex, it is
locally Lipschitz on $D$
and we can calculate its generalized gradient set; see \cite[Lemma 4]{Kon06f}.

%5=========================lm 5.2====================================================
\begin{lemma} \label{lm:5.2} {\em Let the assumptions in {\bf (A1${}''$)} be fulfilled. Then,
at any point $x\in D$, there exists the generalized gradient set
$$
\partial^{\uparrow} \varphi _{\alpha }(x) =
   G(x)-\left[ \partial^{\uparrow} G(x)^{\top}-\alpha I \right] (y_{\alpha }(x)-x),
$$
where $\partial^{\uparrow} G(x)$ denotes the generalized Jacobian of $G$ at $x$.}
\end{lemma}

We recall that a mapping $G :X\rightarrow \mathbb{R}^{n}$ is said to be

(a) {\em monotone} if, for each pair of points $x, y\in X$, it holds
that
$$
\langle G(x)-G(y), x-y\rangle \geq 0;
$$

(b)  {\em strongly monotone} with constant $\tau > 0$ if, for each
pair of points $x, y\in X$, it holds that
$$
\langle G(x)-G(y), x-y\rangle \geq \tau \|x-y\|^{2}.
$$

Let us take the strong monotonicity assumption on the mapping
$G$.

{\bf (A4)} {\em The mapping $G : D \rightarrow \mathbb{R}^{n}$ is
strongly monotone with constant $\tau >0$}.

The strong monotonicity enables us to obtain the desired coercivity and stationarity
properties; see \cite[Lemmas 5 and 6]{Kon06f}.

%5=========================lm 5.3====================================================
\begin{lemma} \label{lm:5.3}  {\em Let the assumptions in {\bf (A1${}''$)} and {\bf (A4)} be fulfilled. Then:

(a)  VI (\ref{eq:5.1}) has a unique solution;

(b)  There exists a number $\sigma >0$ such
that
$$
\varphi _{\alpha }(x)\geq \sigma \|x-x^{*}\|^{2} \quad
\forall x\in D,
$$
where $x^{*}$ is a unique solution to VI (\ref{eq:5.1});

(c) For each point $x\in D$ and for all elements $V \in
\partial^{\uparrow} G(x)$ it holds that
$$
\langle G(x)-(V^{\top}-\alpha I)(y_{\alpha }(x)-x), y_{\alpha }(x)-x \rangle
\leq -\tau \|y_{\alpha }(x)-x\|^{2}.
$$}
\end{lemma}

Therefore, the assumptions in {\bf (A2)} and {\bf (A3)} are fulfilled with $f=\varphi _{\alpha }$,
besides, $D^{*}=D^{0}=D^{e}$, and we can
obtain the convergence result for (GFPMS) directly from Theorem \ref{thm:3.1}.

%3=========================thm 5.1============================================

\begin{theorem} \label{thm:5.1}{\em
Let assumptions {\bf (A1${}''$)} and {\bf (A4)} be fulfilled.
If we apply (GFPMS) with $\beta< \tau$, then the following assertions are true:

 (i) The termination of (GFPMS) yields a unique solution to VI (\ref{eq:5.1}).

(ii) The sequence $\{x^k \}$ converges to a unique solution to VI (\ref{eq:5.1}).}
\end{theorem}

There exist various gap function based methods with
line-search procedures  for different classes of VIs; see e.g. \cite{Pat99,FP03,Kon12,Kon13}.
Again, we can create new versions of these methods
after the proper mapping substitution in (SBM).

%66666666666666666666666666666666666666666666666666666666666666666666666666666666666666666666666666

\section{Computational Experiments}\label{sc:6}

In order to check the performance of the proposed methods we carried
out computational experiments.  The main goal was to compare them with the
methods having the same direction mapping but utilizing
the Armijo line-search. For more clarity, we describe now
the corresponding modification of (SBM).

%===========================Method (SBM)==================================
\medskip
\noindent {\bf Method (ABM).}

{\em Step 0:} Choose a point $x^{0}\in D_{\gamma}$, numbers $\beta \in (0,1)$ and $\theta \in (0,1)$.
Set $k=0$.

{\em Step 1:} Take a point $y^{k}=y(x^k)$.
If $y^{k}=x^{k}$, stop. Otherwise set $d^{k}=y^{k}-x^{k}$.

{\em Step 2:} Determine $m$ as the smallest nonnegative integer such that
$$
 f (x^{k}+\theta ^{m} d^{k}) \leq f (x^{k})-\beta \theta ^{m} \|d^{k}\|^{2},
$$
set $\lambda_{k}=\theta ^{m}$, $x^{k+1}=x^{k}+\lambda_{k}d^{k}$,
$k=k+1$ and go to Step 1.
\medskip

Various implementations of this method were investigated in many works;
e.g. see \cite{Pat98,Pat99,Kon12,Kon13} and the references therein.
Its convergence properties are similar to those of the other known descent
methods with line-search and it does not require a priori information.
The presence of this line-search at each iteration is the main difference
from (SBM).  Also, we took for
comparison the known non-monotone method with the divergent series step-size rule,
which does not require line-search or a priori information;
see e.g. \cite[Chapters V and VII]{Pol87}.

We compared all the methods for different dimensionality.
They were implemented in Delphi with double precision
arithmetic. Namely, we indicate the number of iterations (it) and
the total number of calculations of the goal function value (kf) for attaining the same
accuracy $\varepsilon=0.01$ with respect to the error function
$$
\Delta (x)=\| x-\pi_{D}[x-y(x)]\|.
$$
We took $\theta =0.5$ for  (ABM). For
(SBM), we simply set $\lambda_{k+1} =\lambda_{k}$
if  (\ref{eq:3.2a})  holds, and $\lambda_{k+1} =\sigma \lambda_{k}$
with $\sigma = 0.9$  otherwise.

First we applied the methods to smooth convex optimization problems of form (\ref{eq:2.3}).
More precisely, we chose
\begin{equation} \label{eq:6.1}
f(x)= 0.5 \| Px-q \|^{2},
\end{equation}
the elements of the $m \times n $ matrix $P$ were defined by
$$
p_{ij}= \left\{ {
\begin{array}{ll}
\displaystyle
\sin (i) \cos (j) \quad & \mbox{if} \ i \neq j, \\
\sin (i) \cos (j)+2 \quad & \mbox{if}  \ i=j;
\end{array}
} \right.
$$
and
$$
q_{i}= \sum  \limits_{j=1}^{n}  p_{ij}, \ i=1, \dots, m.
$$
We utilized the mapping $y_{\alpha }(x)=\pi_{D}[x-\alpha^{-1} \nabla f(x)]$ with $\alpha=1$
and (GPMS). Analogously, setting $y(x)=y_{1}(x)$ in (ABM), we obtain
the well known gradient projection method with Armijo line-search.
We call it (GPMA) for brevity. We set $\beta =0.5$ for both the methods.
We also implemented the gradient projection method
with the divergent series step-size rule:
$$
\displaystyle x^{k+1}=\pi_{D}[x^k-\lambda_k \nabla f(x^k)], \quad  \lambda_k =1/(k+1),  \ k=0,1,\ldots;
$$
see e.g. \cite[Section 7.2.2]{Pol87}.
We call it (GPMD) for brevity.

In the first series, we took the
feasible set $D =\mathbb{R}^{n}_{+}$ where
$$
\mathbb{R}^{n}_{+}=\{ x \in \mathbb{R}^{n} \ : \ x_{j} \geq 0, \ j=1, \dots, n \}
$$
and the starting point $x^{0}_{j} = n/2+\sin (j) $ for $j=1, \dots, n$.
The results  are given in Table \ref{tbl:1}.
\begin{table}
\caption{Convex optimization: Test 1
({\em it} is the number of iterations, {\em  kf} is
the number of function calculations)} \label{tbl:1}
\begin{center}
\begin{tabular}{|rr|rr|rr|rr|}
\hline
& {}  & (GPMA) & {} & (GPMS) & {} & (GPMD) & {}  \\
\hline
 $m$ & $n$                & {\em  it} & {\em kf }  & {\em  it} & {\em kf } & {\em  it} & {\em kf }  \\
\hline
  $2$ &  $5$            & 4 & 14  &  12 & 21 &  17 & 18  \\
\hline
   $4$ &  $5$            & 15 & 57 & 30  &  35 & 40  &  41   \\
\hline
  $5$ &   $10$           & 18 & 76 & 28  &  47 & -  &  -  \\
\hline
$25$ & $50$           & 344 & 2683 & 637 & 679 & -  &  - \\
\hline
$50$ & $100$         & 1229 & 12025 & 2633  & 2689 & -  &  - \\
\hline
\end{tabular}
\end{center}
\end{table}
(GPMD) showed very slow convergence when $n>5$.
For the case where $m=5$ and $n=10$, it attained only
the accuracy $0.108$ in $5000$ iterations.
For this reason, we made the further comparison only for (GPMA) and (GPMS).

In the second series, we took the same cost function from (\ref{eq:6.1}),
the feasible set
$$
D =\{ x \in \mathbb{R}^{n} \ : \ -5 \leq x_{j} \leq 5, \ j=1, \dots, n \},
$$
and the starting point $x^{0}_{j} = -5$ for $j=1, \dots, n$.
The results  are given in Table \ref{tbl:2}.
\begin{table}
\caption{Convex optimization: Test 2
({\em it} is the number of iterations, {\em  kf} is
the number of function calculations)} \label{tbl:2}
\begin{center}
\begin{tabular}{|rr|rr|rr|}
\hline
& {}  & (GPMA) & {} & (GPMS) & {}  \\
\hline
 $m$ & $n$                &  {\em  it} & {\em kf }   & {\em  it} & {\em kf }  \\
\hline
  $2$ &  $5$            & 4 & 24  &  14 & 21  \\
\hline
   $4$ &  $5$            & 17 & 65 & 35  &  38  \\
\hline
  $5$ &   $10$           & 19 & 80 & 60  &  66  \\
\hline
$25$ & $50$           & 225 & 1778 & 440 & 463  \\
\hline
$50$ & $100$         & 748 & 7445 & 1624 & 1660 \\
\hline
\end{tabular}
\end{center}
\end{table}

Next we applied the methods to variational inequality problems.
We chose the nonlinear strongly monotone mapping in VI (\ref{eq:5.1})
as follows:
$$
G(x)= Ax+b+\mu C(x), \ A=A'+A'',
$$
the elements of the $n \times n $ matrix $A'$ were defined by
$$
a'_{ij}= \left\{ {
\begin{array}{rl}
\displaystyle
\sin (i) \cos (j)/(i+j) \quad & \mbox{if} \ i<j, \\
\sin (j) \cos (i)/(i+j) \quad & \mbox{if} \ i>j, \\
\sum \limits_{s \neq i} | a'_{is}| +2 \quad & \mbox{if} \ i=j;
\end{array}
} \right.
$$
the elements of the $n \times n $ matrix $A''$ were defined by
$$
a''_{ij}= \left\{ {
\begin{array}{rl}
\displaystyle
\sin (ij)\ln(1+i/j) \quad & \mbox{if} \ i<j, \\
-\sin (ij)\ln(1+j/i) \quad & \mbox{if} \ i>j, \\
0 \quad & \mbox{if} \ i=j;
\end{array}
} \right.
$$
and
$$
b_{i}= -10\sum  \limits_{j=1}^{n}  a_{ij}, \ i=1, \dots, n.
$$
This means that the  matrix $A$ is positive definite and asymmetric.
The parameter $\mu$ was set to be $10$,
the mapping $C(x)$ was chosen to be diagonal with the elements
$$
C_{i}(x)= \arctan( x_{i}-2), \ i=1, \dots, n.
$$

We utilized the mapping $ y_{\alpha }(x) =\pi_{D}[x-\alpha^{-1} G(x)]$ with $\alpha=1$
and (GFPMS). Analogously, setting $y(x)=y_{1}(x)$ in (ABM), we obtain
the well known descent projection method with Armijo line-search.
We call it (GFPMA) for brevity. We set $\beta =0.4$ for both the methods.

In the third series, we took the feasible set
$$
D =\{ x \in \mathbb{R}^{n} \ : \ 1 \leq x_{j} \leq 6, \ j=1, \dots, n \},
$$
and starting point $x^{0}_{j} = 6$ for $j=1, \dots, n$.
The results are given in Table \ref{tbl:3}.
\begin{table}
\caption{Variational inequality: Test 3
({\em it} is the number of iterations, {\em  kf} is
the number of function calculations)} \label{tbl:3}
\begin{center}
\begin{tabular}{|r|rr|rr|}
\hline
 & (GFPMA) & {} & (GFPMS) & {}  \\
\hline
  $n$                &  {\em  it} & {\em kf }   & {\em  it} & {\em kf }   \\
\hline
  $5$            & 4 & 14  &  20 & 26  \\
\hline
  $10$            & 8 & 23 & 21  &  27  \\
\hline
  $20$           & 14 & 48 & 40  &  45  \\
\hline
 $50$           & 47 & 161 & 48 & 53  \\
\hline
 $100$         & 85 & 320 & 92 & 97 \\
\hline
$200$         & 148 & 660 & 145 & 150 \\
\hline
$500$         & 375 & 2143 & 345 & 351 \\
\hline
$1000$         & 761 & 5076 & 708 & 716 \\
\hline
\end{tabular}
\end{center}
\end{table}

In almost all the cases, the implementations of (SBM), which do not use line-search, showed
rather rapid convergence, they outperformed the implementations of
(ABM) in the total number of goal function calculations.

%77777777777777777777777777777777777777777777777777777777777777777777777777

\section{Conclusions}

We suggested a new simple adaptive step-size procedure in a general class of
solution methods for optimization problems, 
whose goal function may be non-smooth and non-convex.
This procedure does not require any line-search or
a priori information, but takes into account behavior of the iteration sequence.
Therefore, it reduces the implementation cost of each iteration essentially
in comparison with the descent line-search methods.
We established convergence of the method under mild assumptions involving the
usual coercivity condition. We showed that this new procedure yields in fact a general
framework for optimization methods. In particular,
a new gradient projection method for smooth constrained optimization problems
and a new projection type method for minimization of the gap function of a
general variational inequality can be obtained within this framework.
The preliminary results of computational tests showed efficiency
of the new procedure.

%%%%%%%%%%%%%%%%%%%%%%%%%%%%%%%%%%%%%%%%%%%%%%%%%%%%%%%%%%%%%%%%%%%%%%%%%%%%%%%%%%%%%

\section*{Acknowledgement}

The results of this work were obtained within the state assignment of the
Ministry of Science and Education of Russia, project No. 1.460.2016/1.4.
In this work, the author was also supported by Russian Foundation for Basic Research, project No.
16-01-00109a.

%##########################################################################

\end{document}